\documentclass[11pt]{amsart} 
\usepackage{amsthm,amsbsy,amsmath,amssymb,amscd,amsfonts,array}
\usepackage{mathrsfs,verbatim,enumerate,enumitem,color,graphicx}
\usepackage[all,cmtip]{xy}
\xyoption{arc} 
\usepackage{hyperref}
\usepackage{arydshln}
\usepackage{multirow}

\newtheorem{thm}{Theorem}[section]
\newtheorem{prop}[thm]{Proposition}     
\newtheorem{lem}[thm]{Lemma}
\newtheorem{cor}[thm]{Corollary}

\theoremstyle{definition}

\newtheorem{question}[thm]{Question}

\newtheorem{rem}[thm]{Remark}

\numberwithin{equation}{section}

\DeclareFontFamily{OT1}{rsfs}{}
\DeclareFontShape{OT1}{rsfs}{n}{it}{<-> rsfs10}{}
\DeclareMathAlphabet{\curly}{OT1}{rsfs}{n}{it}

\newcommand{\ZZ}{\mathbb{Z}} 
\newcommand{\RR}{\mathbb{R}} 

\newcommand{\Mm}{{\mathbf{m}}}

\newcommand{\om}{\overset{\circ}{m}}
\newcommand{\kk}{\kappa}

\newcommand{\be}{\begin{eqnarray*}}
\newcommand{\ee}{\end{eqnarray*}}
\newcommand{\bne}[1]{\begin{eqnarray} \label{#1} }
\newcommand{\ene}{\end{eqnarray}}

\newcommand{\bp}{
   \arraycolsep=6pt 
   \def\arraystretch{1}
   \begin{pmatrix}  
}
\newcommand{\ep}{
   \end{pmatrix}
   \arraycolsep=2pt 
   \def\arraystretch{1.2}
}

\newcommand{\toz}{\operatorname{toz}}

\newcommand{\SL}{\operatorname{SL}}

\def\arraystretch{1.2}

\begin{document}


\author{Fer\.{ı}t \"{O}zt\"{u}rk}
\address{Bo\u{g}az\.{ı}\c{c}\.{ı} \"{U}n\.{ı}vers\.{ı}tes\.{ı}, Department of Mathematics, \.Istanbul, Turkey}
\email{ferit.ozturk@boun.edu.tr}

\title[Curves on surfaces with prescribed  intersections]{Curves on surfaces with prescribed pairwise intersection numbers}

\begin{abstract} 
Given an ordered sequence of $N$~choose~2 integers, we give necessary and sufficient conditions to have an ordered collection of $N$ simple closed curves on a torus such that the algebraic pairwise intersections of those curves are the given integers. 
We also present partial answers towards realizing curve systems on a surface with minimal genus, given the intersection numbers.
\end{abstract}

\maketitle

\section{Introduction} \label{section:intro}

For $N\in\ZZ, N\geq 2$, given an ordered $N \choose 2$-tuple of integers  
$$\Mm=\left( m_{12}; m_{13}, m_{23}; m_{14}, m_{24}, m_{34}; \ldots; m_{1,N}, m_{2,N}, \ldots, m_{(N-1),N}\right)$$ 
we want to find a closed connected oriented 2-manifold (shortly a surface) $\Sigma_g$ of an arbitrary genus $g$ 
and a collection of $N$ oriented, connected, simple closed  curves $\gamma_1,\ldots,\gamma_N$ (called a curve system) on 
$\Sigma_g$ such that the
pairwise algebraic (i.e. signed) intersections of the curves realize the given intersection scheme $\Mm$, i.e. 
$\gamma_i\cdot \gamma_j=m_{ij}$ for all $1\leq i<j \leq N$. (Similar question for geometric (i.e. unsigned) intersections 
in case of nonegative $m_{ij}$'s.) We shortly call such an $\Mm$ an  $N$-scheme.
The question is to find the {\it least} possible such genus for the given $\Mm$, which we denote by $g(\Mm)$.
It can be shown immediately that any scheme is realized on some surface with sufficiently large genus;
i.e. $g(\Mm)$ is finite (see Section~\ref{section:prel}).

A natural question regarding placing curves on surfaces that has been studied extensively is the curve packing problem: given an upper bound for the number of pairwise (geometric) intersections the maximal number of essential, nonhomotopic simple closed curves is to be determined. In \cite{JMM} an upper bound was given for the maximal value, which is linear in the given bound on a torus. This bound on a torus was improved considerably in \cite{Ag1},\cite{Ag2}, with a direct application to knot theory and 3-manifolds.
In \cite{Pr} bounds for maximal number of curves intersecting at most once were given for punctured hyperbolic surfaces.   
Those bounds have been improved for closed surfaces in \cite{Gr}.

There are also results regarding the curve graphs on a torus.
The curve graph on a torus (sometimes called the Farey graph) 
is a graph $\mathcal{F}$ which has a vertex for each isotopy class of curves; 
two distinct vertices are connected by an edge if and only if the corresponding nonoriented classes intersect minimally (i.e. once on a torus). 
More generally given $m\in\ZZ^+$ one can define graphs where two vertices are connected by an edge if and only if 
the absolute value of the algebraic intersection is exactly $m$ or is $\leq m$ respectively; these graphs are called $m$-Farey graphs, 
denoted by $\mathcal{F}_m$ and $\mathcal{F}_{\leq m}$ respectively \cite{GLRRX}. 
Results of curve packing problems on tori can be related directly with properties of Farey graphs. For instance an upper bound for the maximal number of curves on a torus that intersect pairwise at most $m$ times is an upper bound for the chromatic number of the graph $\mathcal{F}_{\leq m}$. Many combinatorial properties of $\mathcal{F}_m$ and $\mathcal{F}_{\leq m}$  have been worked out in the last couple of decades;
see e.g. \cite{ABG}, \cite{AG}, \cite{GLRRX}.

These works are related with the crossing numbers of planar (multi)graphs too. There are various results under the title of crossing lemmas; see e.g. \cite{HP} or \cite{KPTU} and references therein.
 
Instead of going along those directions, in this small note we set off to construct the curve systems directly for a given scheme.
Using simple manipulations we detect the necessary and sufficient conditions to realize a given scheme on a torus. As seen in the sequel, although the ideas are simple, we believe this definitive result has not been announced before. We are persuaded that our clear-cut result below will be of use in the aforementioned collection of problems above.

\begin{thm}
\label{torus}
An $N$-scheme $\Mm$ with all $m_{ij}\neq 0$ is realized on a torus if and only if the following three conditions are satisfied:
\begin{align}
&\label{triangle} \tag{$\triangle$} 
 \mbox{For any index triple $i,j,k$, $(m_{ij},m_{ik})=(m_{ij},m_{jk})=(m_{ik},m_{jk})$.} 
 \\
& \mbox{For any index quadruple $i,j,k,l$, $m_{ij}m_{kl}-m_{ik}m_{jl}+m_{il}m_{jk}=0$.} \label{square} \tag{$\square$}  \\
& \toz(\Mm;p)<p \mbox{ for each prime $p < N$.} \label{circledast} \tag{$\circledast$} 
\end{align}
\end{thm}
The special case when some $m_{ij}=0$ is dealt with in Section~\ref{section:prel}.
We define the rational number $\toz(\Mm;p)$ in Section~\ref{section:toz}. It is an integer assigned to $\Mm$ using the divisibility data of $m_{ij}$'s by the prime $p$, more precisely on the $p$-valuation $\nu_p(\cdot)$ of the entries.  

The proof of Theorem~\ref{torus} is in Section~\ref{section:g1}. The problem takes place in the exterior algebra with two generators. We must basically find  the integer solutions to a linear system since on a torus the isotopy classes of curves are determined by pairs of integers. The condition \ref{square} in Theorem~\ref{torus} is readily recognized as the Grassmann-Plücker relations. Indeed our problem can be viewed as finding an integral basis for a 2-plane in $\RR^N$ where the plane is given via the Plücker coordinates $\Mm$ (in our case they are not homogeneous coordinates). A necessary and sufficient condition for $\Mm$ to describe a plane in $\RR^N$ is the Grassmann-Plücker relations \ref{square}, i.e. $\Mm$ lies on the Plücker quadric. 

The further requirement in the problem is that an integral basis for this plane must have each corresponding pair of coordinates relatively prime. This gives the problem the number theoretic flavor, which introduces the conditions \ref{triangle} and \ref{circledast}.

We will see in the sequel that in special and simple cases \ref{circledast} is satisfied readily. For instance,

\begin{cor} \label{kolayca}
Assume an $N$-scheme $\Mm$ satisfies \ref{triangle} and \ref{square}. Denote the common pairwise greatest common divisor in \ref{triangle} by $g_{ijk}$. Then 
$\Mm$ is realized on a torus if there are three indices $1\leq i < j < k \leq N$ such that one of the following sufficient conditions is satisfied:
\begin{enumerate}
\item $g_{ijk}=1$.
\item $\nu_p(g_{ijk}) = 0$ for every prime $p\leq N$.
\end{enumerate}
In the other direction, $\Mm$ is not realized on a torus if for some prime $p < N$
$\nu_p(m_{ij}) > 0$ is constant regardless of $i,j$.
\end{cor}

Given $N$, let $\mathfrak{g}(N)$ denote the minimum genus which realizes all possible  $N$-schemes, i.e. 
$$ \mathfrak{g}(N) = \sup_{N\textrm{-schemes } \Mm}\!\! g(\Mm).$$
This is a finite integer. Indeed we show in Section~\ref{section:prel} that $ \mathfrak{g}(N)$ is bounded above by $-2+N(N+1)/2$.  
Of course this is a loose bound. For instance $\mathfrak{g}(2)=1$. 
Also it follows from Theorem~\ref{torus} that $\mathfrak{g}(3)>1$, i.e. there are integer triples that cannot be realized by any curve triple on a torus.
However some 3-schemes may be expressed as a sum of two 3-schemes each of which is realizable on a torus, so that $\Mm$ may be realized on a genus-2 surface
obtained by joining the two tori via a connected sum.
If a 3-scheme $\Mm$ cannot be expressed in such a way while $g(\Mm)=2$, we say that $\Mm$ is 2-endemic.
(see Section~\ref{section:g1} for the more general and precise definition). 
We prove
\begin{thm}
\label{triples}
(a) $\mathfrak{g}(3)=2$, i.e. all 3-schemes are realized on a genus-2 surface. Moreover no 3-scheme
is $2$-endemic. \\
(b) There are infinitely many 2-endemic 4-schemes.
\end{thm}
The second claim of part (a) declares simply that any 3-scheme which is realized on $\Sigma_2$
is so because it has two summands each of which is realized on a torus. We prove part (a) in Section~\ref{section:tria}
and we give the construction of 2-endemic 4-schemes for part (b) in Section~\ref{section:endem}.

Before proceeding let us state the apparent question.

\begin{question}
What is the characterization of all schemes which can be realized on a genus-2 surface? 
\end{question}

Of course this is a difficult and involved question, which can be asked for an arbitrary genus$\;>2$ too. However the complications are already there in genus 2. \\

{\bf Acknowledgements:} The core of this work was built during the last weeks of the author's visit to {\em Le Laboratoire de Mathématiques Jean Leray, Nantes Université} and
his first weeks of his visit to  {\em the Erd\H{o}s Center,  Alfréd Rényi Institute of Mathematics} between September 2022--July 2023. He would like to thank  both institutions for the great working environment and hospitality.

\section{Curve systems on a torus} \label{section:g1}

\subsection{Preliminaries.} \label{section:prel}
\begin{figure}[h]
\centering
\includegraphics[width=0.7\columnwidth]{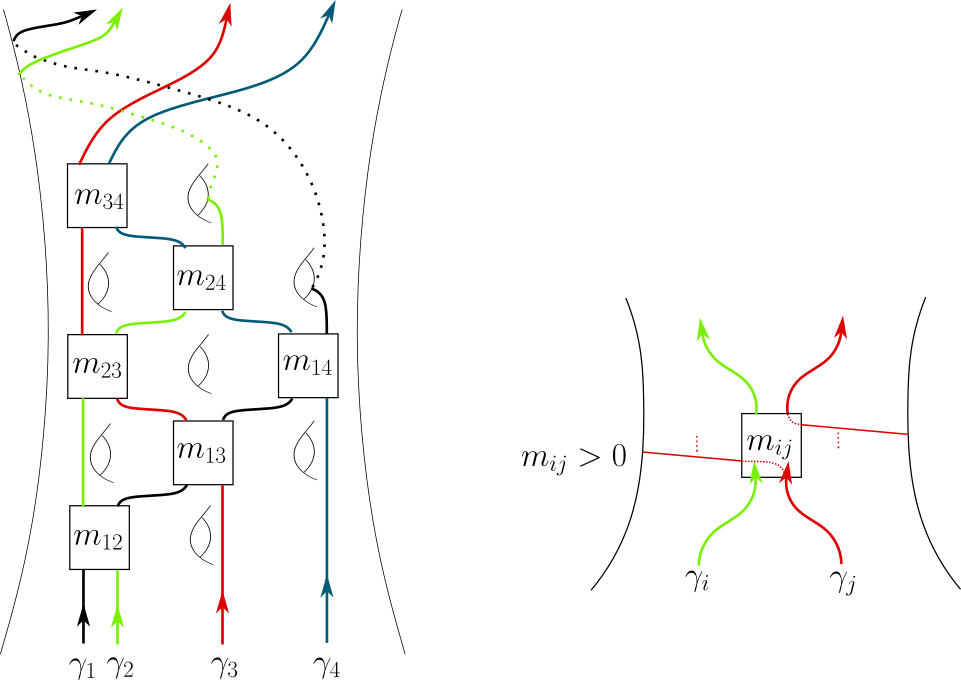}
\caption{Any given 4-scheme is realized on a genus-8 surface. Top and bottom are identified in the left picture. Each intersection $m_{ij}$ is realized on a distinct annulus as depicted on the right hand side for positive $m_{ij}$.}
\label{babayuzey}
\end{figure}

Any $N$-scheme $\Mm$ is realized on some surface with sufficiently large genus;
i.e. $g(\Mm)$ is finite. In fact one can construct a surface as in Figure~\ref{babayuzey}. 
Note that there the curves are isotopically distinct. Moreover since each pair of curves is minimally intersecting,
$m_{ij}$ is the minimum signed intersection possible between the isotopy classes of $\gamma_i$ and $\gamma_j$.
The genus of the surface constructed is $-2+N(N+1)/2$. 
Of course it should be possible to decrease this loose upper bound on genus.

Suppose that the given scheme $\Mm$ can be realized on a torus.
If $m_{ij}=0$ for some $i<j$ then either the corresponding curves are isotopic so that
their pairwise intersection with each other curve is the same; or, say, $\gamma_i$ is empty so that its intersection with any other curve is 0 too.
In both cases we see that the realizability of that $\Mm$ on a torus is equivalent to the realizability of the scheme with $\gamma_i$ ignored.
Because of this particular and basic feature of 0 intersection on a torus, in the following sections we deal with the cases when all $m_{ij}$'s are nonzero.

For example the realization of the scheme $(a;0,0;0,0,0)$, $a$ nonzero, will be a pair of curves on a torus that intersect $a$ times.
The solution for $(a;a,0)$ consists of three curves on a torus with $\gamma_2$ and $\gamma_3$ parallel and
$\gamma_1$ intersects each $a$ times. Note that $(a;b,0)$ with nonzero $a>b$ is not realizable on a torus. 

In general on $\Sigma_g, g>1$,  $m_{ij}=0$ occurs if and only if in the curve complex of $\Sigma_g$ the vertices that represent the
isotopy class of the corresponding curves $\gamma_i,\gamma_j$ are connected by an edge.
For example $\gamma_i,\gamma_j$ may be homologous but not isotopic.

A priori some $N$-schemes are simpler than the others in the following sense.
Suppose that the $N$-schemes $\Mm'$ and $\Mm''$ are realized by the curves $\gamma'_1,\ldots,\gamma'_N$ and $\gamma''_1,\ldots,\gamma''_N$ 
on some surfaces $\Sigma_{g'}$ and $\Sigma_{g''}$ of positive genera $g'$ and $g''$ respectively
and that there are oriented arcs $\omega'\subset \Sigma_{g'}$ and $\omega''\subset \Sigma_{g''}$ such that
$\omega'$ intersects $\gamma'_i$'s in the same order but each with the opposite sign (or each with the same sign) as $\omega''$ intersects $\gamma''_i$'s ($1\leq i\leq N$).
Then the $N$-scheme $\Mm = \Mm' + \Mm''$ is realized on a surface $\Sigma_{g}$ with genus $g=g'+g''$. In fact, 
$\Sigma_{g}$ can be constructed as a
connected sum of $\Sigma_{g'}$ and $\Sigma_{g''}$ at a pair of disks that are sufficiently small neighborhoods of the arcs 
$\omega'$ and $\omega''$. In that way each curve $\gamma'_i$ on $\Sigma_{g'}$ is connected to the corresponding curve on $\Sigma_{g''}$.
(In case $\omega''$ intersects with same signs then the corresponding curves are reversed orientation.)
This final curve $\gamma_i$ is homologous to the disjoint union of the two curves on $\Sigma_{g}$.
Since the algebraic intersection is an invariant on homology classes, it follows that $\Mm$ is realized by $\gamma_i$'s.

For example for the 3-scheme $(a;b,0)$  that we mentioned above, which is not realizable on a torus,
one can express $(a;b,0)$ as $(b;b,0)+(a-b;0,0)$. Here both 3-schemes on the right hand side can be realized on a torus 
and $\gamma'_1,\gamma'_2$ on one torus can be connected to $\gamma''_1,\gamma''_2$ on the other so that
$(a;b,0)$ can be realized on $\Sigma_2$.

The construction in the reverse direction is not valid in general. We will show in Section~\ref{section:endem} that 
there are schemes realized on some $\Sigma_g$ ($g=2$) for which there can be no such  $g',g'',\Mm'$ and $\Mm''$ 
as above. We call such a scheme $g$-endemic.

\subsection{Definition of $\toz$.}\footnote{$\toz$ is a Turkish word which means \it{dust}.}  \label{section:toz}
Given an $N$-scheme $\Mm$, let the indices $1,2,3$ satisfy the condition \ref{triangle}.
We denote by $g_{123}$ the common value $(m_{12},m_{13}) = (m_{12},m_{23}) = (m_{13},m_{23})$. In case $g_{123}>1$, we denote its distinct prime divisors by $g_l, l=1,\ldots,s$, and reserve the letter $l$ for this index in the sequel. For each $g_l$
we denote the $g_l$-valuation $\nu_{g_l}(m_{ij})$ of $m_{ij}$ by $\nu_{lij}$.
For each $j\geq 2$, we set $\nu_{lj}=\min_{1\leq i < \min(3,j)} \nu_{lij}$
and $\nu_l=\nu_{g_l}(g_{123})$.

For purposes to be clear soon, we define the following integers. 
If a prime $p$ does not divide $g_{123}$ then we set $\toz(\Mm;p)=0$. 
For each $4\leq j\leq N$,
we define the rational number
\begin{equation} \label{toz}
\toz(\Mm;g_l,j)= \left\{ \begin{array}{ll} 
g_l^{\nu_{lj}-\nu_{l}}, & 0<\nu_{l1j}=\nu_{l2j}=\nu_{l3j}\leq \nu_{l}, \\
0 & \mbox{otherwise}.
\end{array} \right.  \end{equation}
We define $\toz(\Mm;g_l,3)=1$ if $0<\nu_{l13}=\nu_{l23}=\nu_{l12}$ and 0 otherwise. 
We set $\toz(\Mm;g_l,2)=1$. 

Finally we define
$$ \toz(\Mm;g_l) = \sum_{j=2}^N  \toz(\Mm;g_l,j).
$$
Note that for each $g_l$ and $j$, $0<\toz(\Mm;g_l,j) \leq 1$ so that $\toz(\Mm;g_l) < N$.

Let us understand this definition in simple cases with the assumption \ref{triangle}. Suppose $N=3$. Then for $g_l=2$, $\toz(\Mm;2) = 2$ if and only if $\nu_{l12}=\nu_{l13}=\nu_{l23}$. This condition is exactly \ref{circledast2} in Section~\ref{section:tria} (together with \ref{triangle}) and it is the only situation when the hypothesis \ref{circledast} of Theorem~\ref{torus} fails for $N=3$.

For $N=4$ and $g_l=3$, $\toz(\Mm;3) = 3$ if and only if $\nu_{lij}$ are the same for each $1\leq i<j\leq 4$. This condition is exactly \ref{circledast3} in Section~\ref{tessera} and it  fails the hypothesis \ref{circledast} of Theorem~\ref{torus} for $N=4$ and $g_l=3$.

For a display of examples of $\toz$ computation, see Section~\ref{section:varexamp}. Now we are ready to tell why Corollary~\ref{kolayca} is true.\\

\noindent {\em Proof of Corollary~\ref{kolayca}.} In case of (1) $\toz(\Mm;p)=0$ for every prime $p$. In case of (2) $\toz(\Mm;p)=0$ for every prime $p<N$. In both cases \ref{circledast} holds. In the other direction, if there is a prime $p < N$ for which 
$\nu_p(m_{ij})>0$ and constant regardless of $i, j$ then $\toz(\Mm;p)=N-1\geq p$. 
\hfill $\Box$

\subsection{Two curves.} 
Given $m\neq 0$ we want to 
find $\gamma_1=[p,q]$, $\gamma_2=[r,s]$ such that $p,q,r,s\in\ZZ$ and each pair is relatively prime.
Suppose we have one solution satisfying 
\begin{align} \label{coz} [p,q]\cdot[r,s]=ps-rq=m. \end{align} 
Denote by $S$ the matrix with columns the solutions. 
Then for any $A\in\SL(2,\ZZ)$ we have $A\cdot S$ a solution too. Converse is not true: not every two solution matrices are carried one to the other 
by an element of $\SL(2,\ZZ)$.

We want to produce all solutions for (\ref{coz}). First we observe that 
$(s,q)|m$. Since $(p,q)=1$, we also get $(q,m)|s$; 
and similarly $(s,m)|q$. Hence we obtain:
\begin{equation}
(s,q)=(q,m)=(s,m), \tag{$\triangle_0$} \label{triangle0}
\end{equation}
which is obviously a necessary condition to have a solution. Denote the common value here by $g$.
We can parametrize all solutions for (\ref{coz}) as follows. First set $s'=s/g, q'=q/g$ and $m'=m/g$.
\begin{lem} \label{paramsoln}
Given $m$, assume $s,q$ are given satisfying \ref{triangle0}. Fix a value for $x$ (and hence for $y$) satisfying $x s' - y q' = 1$.
Then $p,r$  solve (\ref{coz}) if and only if they are of the form
\begin{align} \label{pr} 
p= x m' + \kk q', \\
\nonumber r =y m' + \kk s',
\end{align}
with $\kk\in\ZZ$.
\end{lem}
\noindent {\em Proof.}
We prove the only if direction. The converse is straightforward.

We have $p s'\equiv r q' \;\mathrm{mod}\; m'$.
Then  $p \equiv k q' \;\mathrm{mod}\; m'$ where $k=(s')^{-1} r \;\mathrm{mod}\; m'$.
Similarly $r \equiv k_1 s' \;\mathrm{mod}\; m'$  where $k_1=(q')^{-1} p \;\mathrm{mod}\; m'$. 
But $k\equiv k_1  \;\mathrm{mod}\; m'$.
Therefore $p= x m' + \kk q'$ and $r = y_1 m' + \kk_1 s'$, for some $y_1$ and some $\kk,\kk_1$ that are $k$ mod $m'$.  
But $r = y_1 m' + \kk s' - (\kk_1-\kk)s'= (y_1-\alpha s') m' + \kk s'$ for some $\alpha$. Set $y=y_1-\alpha s'$.
Now these $p,r$ satisfy (\ref{coz}) if and only if $x s'- y q' = 1$. 
\hfill $\Box$ \\

Any solution for $p$ is relatively prime with $q'$; similarly for $r$ and $s'$.
We will investigate the relative primeness of $p$ (and $r$) with $g$ in the following section.

Here we take a short cut to list all relatively prime solutions of  (\ref{coz}).


\begin{prop}
\label{n1cozum}
Given $m\neq 0$, the set of relatively prime pairs of solutions for (\ref{coz}) accepts a free $\SL(2,\ZZ)$ action. 
Each orbit has a representative of the form
$S_r=\begin{pmatrix}
1 & r \\
0 & m
\end{pmatrix}, 0<r<m, (r,m)=1$. Two matrices $S_r$ and $S_{r'}$ are in the same orbit if and only if $r\equiv r' \mod m$. Hence there are $\varphi(m)$ orbits of the action, each in 1-1 correspondence with $\SL(2,\ZZ)$. (Here $\varphi$ is the Euler-$\varphi$ function.)
\end{prop}

\subsection{Three curves.} 
\label{section:tria}
 Let us be given $m_{12},m_{13},m_{23}\in\ZZ-\{0\}$.
For the pair of indices $(i,j)=(1,2)$ and $(1,3)$, we fix a representative 
$S_{ij}=\begin{pmatrix} 1 & r_{j} \\ 0 & m_{1j} \end{pmatrix}$
of an orbit in the solution set as in Proposition~\ref{n1cozum}. The remaining requirement reads:
\begin{align}  \label{or3}   r_{2}m_{13} -  r_{3}m_{12}=m_{23}. \end{align}
So the problem is exactly as in the previous section.
We get the necessary condition 
\begin{equation}
(m_{12},m_{13})=(m_{12},m_{23})=(m_{13},m_{23}). \tag{$\triangle$}
\end{equation}
and denote the common value here by $g_{123}$.
We fix $x,y \in\ZZ$ that solve  
\begin{align} \label{xy} 
x m'_{13} - y m'_{12} = 1
\end{align} 
with $m'_{12}=m_{12}/g_{123}$, and $m'_{13}=m_{13}/g_{123}$ as before.
Then by Lemma~\ref{paramsoln} all solutions for $r_2$ and $r_3$ are 
\begin{align} \label{r23} 
r_{2}= x m'_{23} + \kk m'_{12}, \\
\nonumber r_{3}=y m'_{23} + \kk m'_{13},
\end{align}
with $\kk\in\ZZ$. Throughout our work we reserve the letter $\kk$ to parametrize the solutions of the problem.

To satisfy also the conditions $(r_j,m_{1j})=1, j=2,3$, we claim that in every possible case below, $\kk$
can be chosen accordingly, provided that \ref{triangle} and \ref{circledast} are satisfied. 

Assume \ref{triangle}.
We have already observed that 
$(r_2,m'_{12})=1$ and $(r_3,m'_{13})=1$.
This shows in particular that if $g_{123}=1$ then for any choice of $\kk$, 
the corresponding $r_2$ and $r_3$ are solutions to our problem.

Now assume that  $g_{123}>1$ and let its distinct prime divisors be $g_1, \ldots, g_s$. 
Then the following linear system of inequalities is required to hold along with the equations (\ref{r23}):
\begin{align}
\label{r23nonzero}
r_{2}\not\equiv 0 \;\mathrm{mod}\; g_l, & \\
\nonumber r_{3}\not\equiv 0 \;\mathrm{mod}\; g_l, & \qquad l=1,\ldots,s.
\end{align}
We will turn this system into a system of equivalences for $\kk$ by choosing nonzero right hand sides
and then apply the Chinese Remainder Theorem to find all values for $\kk$, giving in turn values for $r_2$ and $r_3$. 
The inequivalances (\ref{r23nonzero}) above forbid some values of $\kk$. The less forbidden values we get the more solutions for $r_2$ and $r_3$ we have.
So let us investigate the allowed/forbidden values for $\kk$ in every case possible.

\begin{itemize}
\item  Let us  first assume $g_l=2$. To make both $r_2,r_3\equiv 1$~mod~2, at least one of $m'_{12}$, $m'_{13}$, $m'_{23}$ must be even, i.e.
\begin{equation}
 \mbox{not all $\nu_2(m_{ij})$ are the same for $1\leq i<j\leq 3$.} \label{circledast2} \tag{$\circledast_2$}
\end{equation}
For otherwise
by (\ref{xy}) $x$ and $y$ would be of different parities, making  one of  $r_2$ and $r_3$ even for any $\kk$.
As noted just after the definition of $\toz$, \ref{circledast2} is equivalent to \ref{circledast} when $N=3$.

Now, there are several cases. If $m'_{23}$ is even then $\kk$ must be odd.
If  $m'_{12}$ is even  then $r_2$ is always odd. 
To make $r_3$ odd too, $\kk$ will be of different parity than $y$.
Similarly if $m'_{13}$ is even  then $\kk$ will be of different parity than $x$. 
Shortly $\kk$ is to satisfy $\kk\not\equiv xym'_{23}$~mod~$2$.

\item  Now assume $g_l > 2$ and $g_l| m'_{12}$. Then 
$r_2\not\equiv 0$~mod~$g_l$ so that we are left with a single inequality: $r_3\not\equiv 0$~mod~$g_l$.
(Below every equivalence is mod~$g_l$.)
For any $v\not\equiv 0$, the equivalence $r_3\equiv v$ can be 
turned into an  equivalence of the  form  $\kk\equiv u$ since $m'_{13}$ has an inverse  mod~$g_l$, namely $x$.
Thus we get $\kk\not\equiv - (m'_{13})^{-1} y m'_{23}  = -xy m'_{23}$~mod~$g_l$.
Similarly for the case $g_l| m'_{13}$, where we get $\kk\not\equiv x y m'_{23}$~mod~$g_l$.

\item  However if $g_l$ divides neither $m'_{12}$ nor $m'_{13}$, for existence of a solution for some $\kk$ at all, 
there must be $u,v\not\equiv 0$ such that the equivalences
$r_2\equiv u\; (g_l)$ and $r_3\equiv v\; (g_l)$ give way to the same value for $\kk$. 
Observe here that if we assume one of $x,y,m'_{23} \equiv 0$ then any $\kk\not\equiv 0$ 
can be taken by a nonzero choice of $u$ or $v$, and the other can be decided accordingly. 
More precisely, if $g_l|m'_{23}$ then any $\kk\not\equiv 0$ is taken by a nonzero choice of $u$ and corresponding $v$.
Else if $g_l|x$ then any $\kk\not\equiv 0,-(m'_{13})^{-1} y m'_{23}$ can be taken.
If $g_l|y$ then $\kk\not\equiv 0,-(m'_{12})^{-1} x m'_{23}$. However note that one can rid these last two cases 
by choosing $x,y$ appropriately in (\ref{xy}).

\item Now suppose none of $x,y,m'_{23},m'_{12},m'_{13}$ is 0  mod~$g_l$.
Solving $r_2\equiv u; \; r_3\equiv v$ for $\kk$ we get $\kk\equiv (x m'_{23}-u)(m'_{12})^{-1}\equiv (y m'_{23}-v)(m'_{13})^{-1}$.
Equivalently $(x m'_{13} - y m'_{12})m'_{23}=m'_{23}\equiv u m'_{13} - v m'_{12}$.
Of course there are infinitely many solutions for the pair $u,v$. Note that 
$u\equiv 0$ if and only if $v\equiv -(m'_{12})^{-1} m'_{23}$; similarly
$v\equiv 0$ if and only if $u\equiv (m'_{13})^{-1} m'_{23}$.
Since $g_l>2$ we can avoid these two cases in choosing $u$ and $v$. So any compatible, admissible choice of $u$ and $v$ gives
a single value for $\kk$ and moreover any value mod~$g_l$ for $\kk$ (including 0) can be produced in this way,
except the two distinct values: $-(m'_{12})^{-1} x m'_{23}$ and $- (m'_{13})^{-1} y m'_{23}$.
\end{itemize}

Let us summarize. For each case above there is an admissible choice for the right hand sides
that leads to a system of equalities $\kk\equiv u_l$~mod~$g_l$.
Thus by the Chinese Remainder Theorem there is a unique solution for $\kk$~mod~$g_{123}$.

Accordingly we obtain all possible solutions for $\kk$ in (\ref{r23}) and (\ref{r23nonzero}).
These give all possible solutions for $r_2$ and $r_3$.
This proves Theorem~\ref{torus} in the case $N=3$. More precisely we have the following

\begin{prop}
\label{herdurum}
(i) A 3-scheme $\Mm$ with $m_{ij}$'s nonzero is realized on a torus if and only if it satisfies \ref{triangle} and \ref{circledast2}. \\
(ii) In that case the set of solutions consists of the simple closed curve triples  
$(1,0), (r_2,m_{12}),(r_3,m_{13})$ and their images 
under $\SL(2,\ZZ)$ action. Here 
$$r_2\equiv(m'_{13})^{-1} m'_{23}\!\!\mod m'_{12}, \; 
r_3\equiv-(m'_{12})^{-1} m'_{23}\!\!\mod m'_{13}.$$ 
(iii) The values $r_2$ mod~$m_{12}$ and $r_3$ mod~$m_{13}$ are determined by the allowed values for $\kk$ in (\ref{r23}).

In case $g_{123}=1$,  any $\kk\in\ZZ$ is allowed.

Otherwise let $g_l$, $l=1,\ldots,s,$ be the distinct prime divisors of $g_{123}$.
Then the allowed values for $\kk$ are the unique solutions~mod~$g_l$ of the system
\begin{align}
\label{chinese}
\kk\equiv u_l~mod~g_l, \;\; l=1,\ldots,s,
\end{align}
where $u_l$ can be chosen of any value except the values 
$A=-(m'_{12})^{-1} x m'_{23}$~mod~$g_l$ and $B=-(m'_{13})^{-1} y m'_{23}$~mod~$g_l$,
whenever the inverses are defined.
More precisely:
\begin{enumerate}
\item If $g_l|m'_{12}$ then $u_l\not\equiv - x y m'_{23}$~mod~$g_l$.
\item If $g_l|m'_{13}$ then $u_l\not\equiv + x y m'_{23}$~mod~$g_l$.
\item If $g_l|m'_{23}$ then $u_l\not\equiv 0$~mod~$g_l$.
\item Otherwise $u_l\not\equiv A,B$~mod~$g_l$. 
\end{enumerate}

\end{prop}

\begin{cor}
\label{forbiddenno}
In the case $g_{123}=1$ there is no forbidden value for $\kk$.
In the case $g_{123}>1$ if a prime divisor $g_l$ of $g_{123}$ divides $m'_{12}m'_{13}m'_{23}$, then there is a single forbidden value mod~$g_l$.
Otherwise there are 2 forbidden values mod~$g_l$.
\end{cor}

The proposition also shows that the conditions \ref{triangle} and 
\ref{circledast2} applied to any index triple are necessary for any $N$-scheme to be realized on a torus.

As an example which does not satisy \ref{circledast2} consider $\Mm=(6,10;14)$. 
We get $5r_2+3r_3=7$ and $r_2=2+3\kk, r_3=-1-5\kk$. If $\kk$ is even (respectively odd) then $(r_2,m_{12})=2$ 
(respectively $(r_3,m_{13})=2$). So there is no solution for $r_2,r_3$ and therefore $g(\Mm)>1$. 
However one can express $\Mm$ as $(1,5;14)+(5,5;0)$. Here each term is realizable on a torus and hence $g(\Mm)=2$.
This argument works in general and is essentially the core of the proof of  Theorem~\ref{triples}(a), which we are ready to present now. \\

\noindent {\em Proof of Theorem~\ref{triples}(a).} Given $\Mm=(m_{12};m_{13},m_{23})$
if any of the entries is zero then, as we discussed in Section~\ref{section:prel}, it can be realized on $\Sigma_g$ with $g=1$ or 2,
where $g=2$ if and only if $\Mm=(a;b,0), a\neq b$. 
So suppose no entry is zero and $\Mm$ does not satisfy \ref{triangle} or \ref{circledast2}.
To prove the theorem it suffices to show that $\Mm$ can be expressed as a sum
$\underline{\mathbf{p}} + \underline{\mathbf{q}}$ where both $\underline{\mathbf{p}}$ and $\underline{\mathbf{q}}$  are realizable on a torus.
Suppose that, after an index change if necessary,  $m_{12}$ and $m_{13}$ are of the same parity. 
Thanks to the following lemma, we know that there is a $\kk\in\ZZ$ for which 
the triple $\underline{\mathbf{p}}=(m_{12}-\kk; m_{13}-\kk, m_{23})$ of integers are pairwise relatively prime.
So $\underline{\mathbf{p}}$ satisfies \ref{triangle} and \ref{circledast2} and therefore is realizable on a torus. 
The triple  $\underline{\mathbf{q}}=(\kk;\kk,0)$ is realizable on a torus too. 
\begin{figure}[h]
\centering
\includegraphics[width=0.6\columnwidth]{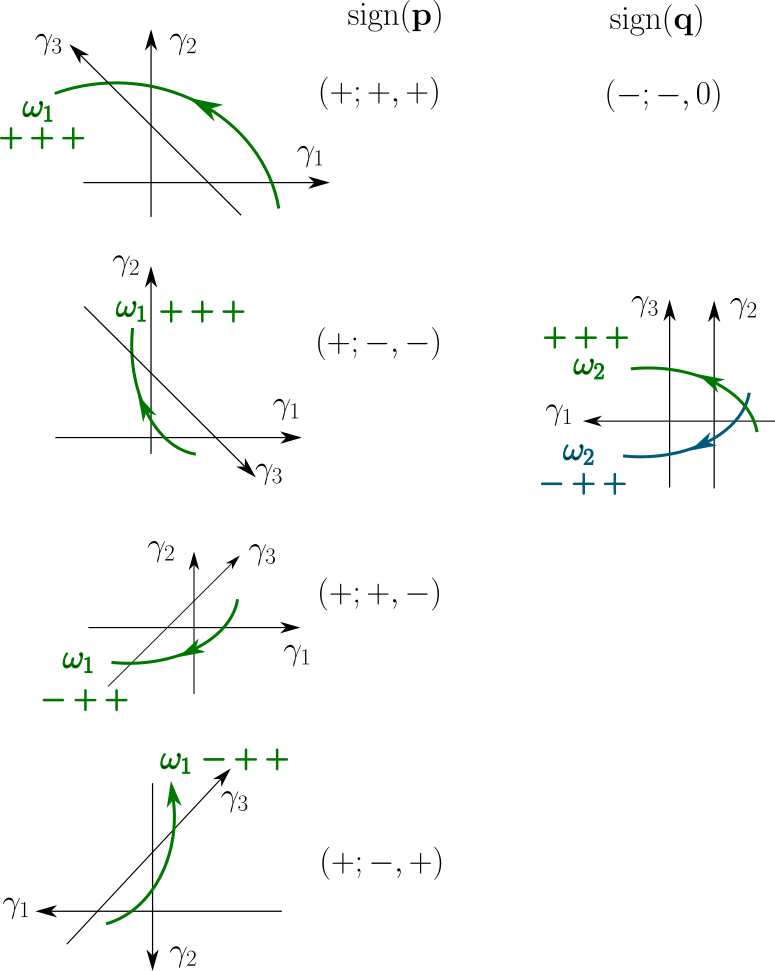}
\caption{For all possible signs of intersections of three curves on $\Sigma_1$ and for the scheme $(\kk;\kk,0)$ ($\kk<0$) on $\Sigma_2$, 
one can find two arcs $\omega_1$ and  $\omega_2$ with appropriate intersection order and signs, as noted in green.}
\label{3connected}
\end{figure}

Moreover choosing a, say, negative $\kk$ one can choose arcs $\omega_1$ and $\omega_2$ on these tori such that they intersect the curves in the same order but
with opposite signs as in Figure~\ref{3connected}. Hence $\Mm$ can be realized on the connected sum of these two tori. 
\hfill $\Box$

\begin{lem}
Let $a,b,c\in\ZZ$ with $a\equiv b \mod 2$.
Then there is $\kk\in\ZZ$ for which the integers $a-\kk,b-\kk,c$  are pairwise relatively prime.
The value for $\kk$ is unique mod the product of distinct prime divisors of $c(a-b)$.
\end{lem}

\noindent {\em Proof.}
To have $(a-\kk,c)=1$, we must have $a-\kk\not\equiv 0 \; (c_j)$ for each prime divisor $c_j$ of $c$; i.e. $\kk\not\equiv a \;  (c_j)$.
Similarly  $\kk\not\equiv b \;  (c_j)$ for each $j$. Moreover we require $(a-\kk,b-\kk)=1$. To get this we must have
$a-\kk\not\equiv 0 \; (d_i)$ for each prime divisor $d_i$ of $a-b$; i.e. $\kk\not\equiv a\; (d_i)$ for each $i$. 
Since $a$ and $b$ are of the same parity, all of the inequalities for $\kk$ for mod~2 are identical. 
So we may create equalities~mod~$c_j,d_i$ for $\kk$ to get a system of linear equations
and apply the Chinese Remainder Theorem to obtain a solution for $\kk$.
\hfill $\Box$

\subsection{Four curves.} 
\label{tessera}
Given a 4-scheme $\Mm$ with no entry zero, assume that the conditions \ref{triangle} and \ref{circledast2} are satisfied for any index triple, 
that $\gamma_1=(1,0)$, and that  solutions for $\gamma_j$'s ($4\geq j>1$) are of the form $(r_{j},m_{1j})$. 
The $r_{j}$'s must satisfy the $\binom{4-1}{2}=3$ equations:
\begin{align*}  
\begin{pmatrix} m_{13} & -m_{12} & 0 \\ 0 & m_{14} & -m_{13} \\ m_{14} & 0 & -m_{12}  \\ \end{pmatrix}
\begin{pmatrix} r_{2} \\ r_{3} \\ r_{4}\end{pmatrix}
= \begin{pmatrix} m_{23} \\ m_{34} \\ m_{24}\end{pmatrix}.
\end{align*}

In order for this system to have a solution, the system must be consistent. Thus we get the necessary condition:
\begin{equation}
\mu_{1234} = m_{12}m_{34}-m_{13}m_{24}+m_{14}m_{23}=0. \label{square4} \tag{$\square_4$}
\end{equation}

\begin{rem} 
From the necessary condition \ref{square4} we may deduce that  all $m_{ij}$ cannot have the same nonzero value, say, $k$.
Moreover in case $k$ is even, (\ref{or3}) has the solutions $(1,0),(r,k),(r+1,k)$; however one of the last two is not a simple curve.
Hence  (\ref{or3}) has no solution.  All this discussion reproves that
if there is a complete subgraph of a $k$-Farey graph $\mathcal{F}_k$ then the number of vertices of that subgraph is at most 3 for $k$ odd, and 2 for $k$ even (cf. \cite[Lemma~3]{GLRRX}).
\end{rem}

Let us now solve the above simplified system for integers $r_2,r_3,r_4$.
As we noted for the case $N=3$, all solutions for $r_2$ and $r_3$ are as in (\ref{r23}).
Similarly we find 
\begin{align} \label{r4} 
m_{12} r_{4} & = r_2 m_{14} - m_{24} & \\
\nonumber  & =  m_{14} (x m'_{23} + \kk m'_{12}) - m_{24} & \textrm{by } (\ref{r23})\\
\nonumber  & =  x (m'_{13} m_{24} - m'_{12} m_{34} ) + \kk m_{14} m'_{12} - m_{24} & \textrm{by } (\square_4)\\
\nonumber  & = ( x m'_{13} - 1) m_{24} + m'_{12} (\kk m_{14}  - x m_{34} ) & \\
\nonumber  & = m'_{12} ( y m_{24} - x m_{34} + \kk m_{14}) & \textrm{by } (\ref{xy}) \\
\label{r4kk} \Longrightarrow g_{123} r_4 & = y m_{24} - x m_{34} + \kk m_{14}.
\end{align}

Thus $r_4$ has an integral solution if and only if there is some $\kk\in\ZZ$ such that 
$g_{123}$ divides $y m_{24} - x m_{34} + \kk m_{14}$.
To elaborate this first we fix a prime divisor $g_l$ of $g_{123}$. We simply observe the following:
\begin{lem}
\label{alldiv}
Assume \ref{triangle}. Below $1\leq i<j$. \\
(a) If $\nu_{li4}>0$ for some $1\leq i\leq 3$, then $\nu_{l4}>0$.  \\
(b) If $\nu_{li4} > \nu_{l}$ for some $1\leq i\leq 3$, then $\nu_{l4} = \nu_{l}$.
Otherwise (i.e. if for every $i$, $\nu_{li4} \leq \nu_{l}$) then $\nu_{l14}=\nu_{l24}=\nu_{l34}=\nu_{l4}$. More precisely up to a permutation of the first three indices
there are 4 possibilities for $\nu_{lij}$'s:
\begin{align*} 
 (v;\nu_{l},\nu_{l};u,v,\nu_{l}), (u;\nu_{l},\nu_{l};v,v,\nu_{l}), (\nu_{l};\nu_{l},\nu_{l};u,\nu_{l},\nu_{l}), (v;\nu_{l},\nu_{l};w,w,w)
 \end{align*}
where $u>v\geq \nu_{l} \geq w$.
\end{lem}
\noindent {\em Proof:} 
Follows immediately from \ref{triangle} through a simple case-by-case analysis.
\hfill $\Box$ \\

Now if $g_{123}=1$ then (\ref{r4}) gives a solution for $r_4$, regardless of the values for $\kk$.

So we assume $g_{123}>1$ and fix a prime divisor $g_l$ of $g_{123}$. 
If $\nu_{l14} > \nu_{l}$ then by the lemma above $\nu_{l4}=\nu_{l}$ so that through (\ref{r4}), $r_4 = - (\om_{12})^{-1} \om_{24}$. 
(Here and below $\om_{ij}=m_{ij}/g_l^{\nu_{lj}}$.)
This value for $r_4$ is independent from  $\kk$, and is nonzero~mod~$g_l$.

When $\nu_{l14} = \nu_{l}$, there are several configurations that arise, as seen in the lemma above. Let us postpone the case when all $\nu_{lij}$ are the same for $j=3,4$. Then we have, for instance, the case $(u;\nu_{l},\nu_{l};\nu_{l},\nu_{l},v)$
with $u,v>\nu_{l}$. Here by (\ref{r4kk}), $r_4=(g_{123}/g_l^{\nu_{l}})^{-1} (y \om_{24} - \kk \om_{14})$~mod~$g_l$. To have $r_4$ nonzero~mod~$g_l$, we must impose
\begin{align}
\label{k-yasak1}  \kk\not\equiv \om_{14}^{-1} y \om_{24} \mod g_l.
\end{align} 
Similar argument is valid for other instances.

Lastly if $\nu_{l14} < \nu_{l}$, then by Lemma~\ref{alldiv} we have the case 
$(u;\nu_{l},\nu_{l};\nu_{l4},\nu_{l4},\nu_{l4})$ with $u\geq\nu_{l}$, up to permuting first three indices. In this case 
for $r_4$ to have a solution in (\ref{r4kk}) we must have
\begin{align}
\label{k-val}  \kk \equiv (\om_{14})^{-1}_{\nu_{l}-\nu_{l4}} (x \om_{34}-y\om_{24})\mod g_l^{\nu_{l}-\nu_{l4}}.
\end{align}
(Here $(b)^{-1}_p$ denotes the inverse of $b$ in $\mathbb{F}_{g_l^p}$.)
Moreover if $\nu_{l4}>0$ we must avoid $r_4\equiv 0 \mod g_l$, so
\begin{align}
\label{k-yasak2} 
\kk\not\equiv (\om_{14})^{-1}_{\nu_{l}-\nu_{l4}+1} (x\om_{34}-y\om_{24}) \mod g_l^{\nu_{l}-\nu_{l4}+1}.
\end{align}
Now we compare these conditions on $\kk$ with the prescribed values for $\kk$ in Proposition~\ref{herdurum}. For instance when $g_l\nmid m'_{13}$, consider the value $B$ 
in Proposition~\ref{herdurum} and its difference with the right hand side of 
(\ref{k-yasak1}) or (\ref{k-val}) or (\ref{k-yasak2}):
\begin{align}
\label{k-ok} 
& (\om_{14})^{-1} (x \om_{34} - y\om_{24}) -(-((m'_{13})^{-1} y m'_{23})) & \\
\nonumber
& \equiv (\om_{14})^{-1} (m'_{13})^{-1} (m'_{13} x \om_{34} - y m'_{13} \om_{24} + y m'_{23} \om_{14}) & \\
\nonumber
& \equiv (\om_{14})^{-1} (m'_{13})^{-1} (m'_{13} x \om_{34} - y m'_{12} \om_{34}) 
& \textrm{by } (\square_4) \\
\nonumber  & \equiv (\om_{14})^{-1} (m'_{13})^{-1} \om_{34} \mod g_l &  \textrm{by } (\ref{xy})
\end{align}
Similar computation runs for the forbidden value $A$ in Proposition~\ref{herdurum}.
In the case of (\ref{k-yasak1}) $\om_{34}\equiv 0 \mod g_l$ so the two forbidden values
(that in  (\ref{k-yasak1}) and $B$ in Proposition~\ref{herdurum}) for $\kk$ are the same. Otherwise in the case of (\ref{k-val}) and (\ref{k-yasak2}) 
$\om_{34}$ is nonzero~mod~$g_l$. Hence in that case (i.e. when $\nu_{l4}<\nu_{l}$), the required value for $\kk$ in (\ref{k-val}) does not coincide with the forbidden values in 
Proposition~\ref{herdurum}. 

The only case we have discarded up to now is the case when $\nu_{lij}$ is a positive constant for $j=3,4$. In that case the forbidden value for $\kk$ in (\ref{k-yasak2}) is different~mod~$g_l$ from those in Proposition~\ref{herdurum}.

Thus for $r_4$ to have an integral solution such that $(r_4,g_{123})=1$ it suffices to
have an admissible value for $\kk$~mod~$g_l$ for each $g_l$. In other words 
the total count of forbidden values~mod~$g_l$ must be less than $g_l$. Taking into account the discussion preceeding (\ref{k-yasak2}) and Corollary~\ref{forbiddenno}, that count would be
$$ \left\{ \begin{array}{ll} 
1, & \nu_{g_l}(m'_{12}m'_{13}m'_{23}) > 0 \mbox{ and } \nu_{li4}\neq \nu_{l}, \exists i=1,2,3; \\ 
2, & \nu_{g_l}(m'_{12}m'_{13}m'_{23}) > 0 \mbox{ and } \nu_{li4}=\nu_{l}, \forall i=1,2,3;\\
2, & \nu_{g_l}(m'_{12}m'_{13}m'_{23})=0 \mbox{ and } \nu_{li4}\neq \nu_{l}, \exists i=1,2,3; \\ 
3, & \nu_{g_l}(m'_{12}m'_{13}m'_{23})=0 \mbox{ and } \nu_{li4}=\nu_{l}, \forall i=1,2,3; \mbox{ i.e. $\nu_{lij}$ is constant.}
\end{array}
\right.
$$
Since we have already assumed that \ref{circledast2} is satisfied for each index triple, 
the total count of forbidden values~mod~$2$ is 1. Meanwhile for $g_l>3$, that count for values~mod~$g_l$ is obviously less than $g_l$. However the count for values~mod~3 can very well be 3, leading no solution. To avoid this we assume that
\begin{equation}
 \mbox{not all $\nu_3(m_{ij})$ are the same for $1\leq i<j\leq 4$.} \label{circledast3} \tag{$\circledast_3$}
\end{equation}

Recall that we have already noted just after the definition of $\toz$ that \ref{circledast3} together with \ref{circledast2} are equivalent to \ref{circledast} when $N=4$.

To finish we still have to show that a prime divisor $p$ of $m_{14}$ that does not divide $g_{123}$ does not divide $r_4$ either. In fact if $p|r_4$ too then 
considering the first line of (\ref{r4}), $p|m_{24}$ hence $p|m_{12}$. 
Now the last line of (\ref{r4}) shows that $p|xm_{34}$. Since $(x,m'_{12})=1$ we
get  $p|m_{34}$ so that $p|m_{13}$ too. This contradicts with $p \nmid g_{123}$.\\

\noindent {\em Proof of Theorem~\ref{torus} for $N=4$.} 
Proposition~\ref{herdurum} shows that \ref{circledast2} is necessary to have a solution.
Moreover it is sufficient if some $m_{ij}$'s are zero. So assume none of them is zero. 
As in Section~\ref{section:tria} we find all solutions for the parameter $\kk$ in (\ref{r23}) 
to get all solutions for $r_2$ and $r_3$ that make the curve triple 
$(1,0),(r_2,m_{12}),(r_3,m_{13})$ solve the problem for $N=3$.

Then we solve for $r_4$ during which \ref{square4} turns out to be necessary. 
As we observed just above, assuming \ref{circledast3} keeps the forbidden values~mod~3 
for $\kk$ less than 3. Now replacing the equivalence $\kk\equiv u_l\mod g_l$ in
(\ref{chinese}) by  the sharper requirements above for each $g_l$,
we employ the Chinese Remainder Theorem to this more restricted system and 
get a solution for $\kk$ and for $r_4$ with $(r_4,m_{14})=1$. 
Hence \ref{triangle}, \ref{circledast2}, \ref{circledast3} and \ref{square4} are also 
sufficient to have a solution.
\hfill $\Box$

\subsection{More curves.}  
It is easy to show inductively that the following situation is true for any $N\geq 5$. 
For $N=5$ proceeding as before, we get  the linear system in $5-1=4$ unknowns and $\binom{4}{2}=6$ equations:
\begin{align*}  
\begin{pmatrix} m_{13} & -m_{12} & 0 & 0 \\ m_{14} & 0 & -m_{12} & 0 \\ m_{15} & 0 & 0 & -m_{12} \\ 
0 & m_{14} & -m_{13} & 0 \\ 0 & m_{15} & 0 & -m_{13} \\ 0 & 0 & m_{15} & -m_{14} \\ 
\end{pmatrix}
\begin{pmatrix} r_{2} \\ r_{3} \\ r_{4} \\ r_{5}\end{pmatrix}
= \begin{pmatrix} m_{23} \\ m_{24} \\ m_{25} \\ m_{34} \\ m_{35} \\ m_{45}\end{pmatrix}.
\end{align*}
Repeating the row operations for the indices 1234 and 1235 which led us to \ref{square4},
we get the independent necessary conditions:
\begin{align*}  
\mu_{1234} =0, \mu_{1235} =0, \mu_{1345} =0.
\end{align*}
Under these conditions, the system above simplifies
to $N-2$ independent equations of the form $m_{1j}r_2-m_{12}r_j=m_{2j}$ and we have 
$\binom{N}{2}- (N-2)$ consistency conditions of the form $\mu_{abcd}=0$. This gives us the necessary condition \ref{square}.

With a more careful analysis we observe that  
\begin{align*}
m_{ae}\mu_{abcd} - m_{ad}\mu_{abce}+m_{ac}\mu_{abde} - m_{ab}\mu_{acde}=0.
\end{align*}
For instance given 5 indices $a,b,c,d,e$, requiring vanishing of only three $\mu$'s, e.g. 
$\mu_{abcd},\mu_{abce},\mu_{abde}$, is equivalent to requiring vanishing of all $\mu$'s with indices among  $a,b,c,d,e$. Therefore  \ref{square} is equivalent to requiring:
\begin{equation}
\mbox{for every pair of indices $1<i,j\leq N$ with $i+1<j$, $\mu_{1,i,i+1,j}=0$.} \label{bsquare} \tag{$\boxtimes$}
\end{equation}
Note that there are $(N-3)(N-2)/2$ such expressions. 

Now assuming \ref{circledast2} and \ref{square} (or equivalently \ref{bsquare}), suppose we solve $r_2$ and $r_3$ as in Section~\ref{section:tria}. 
Then we solve $r_4$ up to  $r_N$ one by one as in Section~\ref{tessera}. 
Here one has to make sure that while solving, say, $r_4$ and $r_5$, the conditions imposed on $\kk$, if any, do not contradict.
In fact, using the notation of the previous section and Lemma~\ref{alldiv}, let us 
fix a prime divisor $g_l$ of $g_{123}$ and consider the difference of the critical values 
for $\kk$ expressed on the right hand side of (\ref{k-val}), or of (\ref{k-yasak2}) whenever they appear:
\begin{align} \label{klerayni}
\alpha-\beta=& \; \om^{-1}_{14} (x \om_{34}-y\om_{24})
- \om^{-1}_{15} (x \om_{35}-y \om_{25}) \\
\nonumber   = & \; \om^{-1}_{14} \om^{-1}_{15} 
\left(x (\om_{15} \om_{34} - \om_{14} \om_{35}) + y (\om_{14} \om_{25}-\om_{15} \om_{24})\right) \\ 
\nonumber  = & \; \om^{-1}_{14} \om^{-1}_{15} 
(x m_{13} {m}_{45} - y m_{12} {m}_{45}) / g_l^{\nu_{l4}+\nu_{l5}}  \\ 
\nonumber  = & \; \om^{-1}_{14} \om^{-1}_{15} {m}_{45} g_{123} / g_l^{\nu_{l4}+\nu_{l5}} 
\end{align}
This manipulation is essentially the same as that in (\ref{k-ok}).

In solving $r_4$ (respectively $r_5$), (\ref{k-val}) and (\ref{k-yasak2}) occur when $\nu_{l}\geq \nu_{l4}=\nu_{li4}, i=1,2,3$ (respectively $\nu_{l}\geq \nu_{l5}=\nu_{li5}$). Without loss of generality we assume $\nu_{l4}\geq \nu_{l5}$.
Note that the  last expression in \ref{klerayni} is an integer; indeed recall that $\nu_{l45}\geq \min(\nu_{l4},\nu_{l5})$ so that $\nu_{l4}+\nu_{l5} \leq \nu_{l}+\nu_{l45}$. 

Now we compare the required/forbidden values for $\kk$ and deduce shortly that no contradiction arises.
First assume that $\nu_{l4}\gneq \nu_{l5}$; this implies $\nu_{45}=\nu_5$. 
In that case $\nu_{g_l}(\alpha-\beta)=\nu_{l}+\nu_{l45} - \nu_{l4} - \nu_{l5}= \nu_{l}-\nu_{l4}$. There can be a contradiction between $\kk\equiv \alpha (g_l^a)$ and
$\kk\not\equiv \beta (g_l^b)$ if and only if $\alpha=\beta$ and $b\leq a$. This latter
condition holds in a single case: when $a=\nu_{l}-\nu_{l5}$ and $b=\nu_{l}-\nu_{l4}+1$.
However then both $a,b>\nu_{l}-\nu_{l4}=\nu_{g_l}(\alpha-\beta)$. Hence $\alpha\not\equiv \beta \mod g_l^a \mbox{ or mod }  g_l^b$; i.e. the required and forbidden values do not contradict. 

Furthermore there can be a contradiction between two required values 
$\kk\equiv \alpha (g_l^a)$ and $\kk\equiv \beta (g_l^b)$. This is the situation   
$\alpha\not\equiv \beta \mod g_l^{\min(a,b)}$. In our case $\min(a,b)=\nu_{l}-\nu_{l4}$.
But $\nu_{l}-\nu_{l4}=\nu_{g_l}(\alpha-\beta)$. So in this case $\alpha\equiv \beta \mod g_l^a$
and no contradiction arises.

The remaining case is when $\nu_{l4} = \nu_{l5}$, so that $\nu_{45}\geq \nu_4$.
In this case $\nu_{g_l}(\alpha-\beta)=\nu_{l}+\nu_{l45} - 2\nu_{l4}$. 
Let us compare the requirements $\kk\equiv \alpha (g_l^a)$ and
$\kk\equiv \beta (g_l^b)$. They occur when $a=b=\nu_{l} - \nu_{l4}$.
A contradiction arises if and only if $\alpha\not\equiv\beta \mod g_l^{\nu_{l} - \nu_{l4}}$,
i.e. when $\nu_{l} - \nu_{l4}>\nu_{g_l}(\alpha-\beta)=\nu_{l}+\nu_{l45} - 2\nu_{l4}
\Leftrightarrow \nu_{l4}>\nu_{l45}$, which is impossible.

Similarly we compare the relations $\kk\equiv \alpha (g_l^a)$ and
$\kk\not\equiv \beta (g_l^b)$. In this case $a=\nu_{l} - \nu_{l4}=b-1$.
A contradiction arises if and only if $\alpha\equiv\beta \mod g_l^a$ and $b\leq a$.
However this latter condition does not hold in this case. 
 
Now let us compare and count the forbidden values in $\kk\not\equiv \alpha (g_l^a)$ and
$\kk\not\equiv \beta (g_l^b)$. The forbidden values are distinct if and only if $a=b$
and $\alpha\not\equiv\beta (g_l^{a})$. Here $a=\nu_3-\nu_4+1$  and $b=\nu_3-\nu_5+1$ so that $a=b\Leftrightarrow \nu_4=\nu_5$ and $\nu_{l45} \geq  \nu_{l4}$. Then $\alpha\not\equiv\beta \mod g_l^a
\Leftrightarrow \nu_{g_l}(\alpha-\beta)<a
\Leftrightarrow \nu_{l}+\nu_{l45} - 2\nu_{l4} < \nu_{l}- \nu_{l4} +1
\Leftrightarrow \nu_{l45} \leq  \nu_{l4}$. Hence 
the forbidden values are distinct mod~$g_l^{\nu_{l}-\nu_{l4}+1}$ if and only if $\nu_{l45}=\nu_{l4}$.

Let us recapitulate: A forbidden value occurs for index $j$ if and only if  $\nu_{lij}=\nu_{lj}, 1\leq i\leq 3$. It is a single value mod~$g_l^{\nu_{l}-\nu_{lj}+1}$.  
The forbidden values corresponding to indices $j$ and $k$ are distinct if and only if $\nu_{lij}=\nu_{lik}=\nu_{ljk}, 1\leq i\leq 3$. They are distinct mod~$g_l^{\nu_{l}-\nu_{lj}+1}$.\\

\noindent {\em Proof of Theorem~\ref{torus}}
We solve $r_j$ for each $4 \leq j\leq N$ as in the previous section, assuming the necessary conditions \ref{circledast2} (for each index triple $1,2,j$)
and \ref{circledast3} (for each index quadruple $1,2,3,j$) and \ref{square}.
This possibly imposes conditions on $\kk$ modulo some powers of each prime divisor $g_l$ of $g_{123}$. As we have discussed just above, any required value for $\kk$ modulo some power of $g_l$ does not contradict with any other required or forbidden value. Furthermore in solving each $r_j$ a distinct forbidden value for $\kk$ is introduced mod~$g_l^{\nu_{l}-\nu_{lj}+1}$ if and only if 
$\nu_{l}\geq \nu_{lj}=\nu_{l1j}=\nu_{l2j}=\nu_{l3j}(=\nu_{l4j})$.
There is an allowed value for $\kk$ if and only if $\toz(\Mm;g_l)<g_l$. 
Hence \ref{circledast} together with \ref{triangle} and \ref{square} are also sufficient to have a solution.
\hfill $\Box$

Let us emphasize that this process above produces all possible solutions to the realization problem, with the first curve fixed as $(1,0)$. The complete solution set is the set of curve systems obtained under the action of $SL(2,\ZZ)$. The fact we produce all possible solutions leads us to the following curious consequence of combinatorial flavor:

\begin{cor} Consider an $N$-scheme $\Mm$ with nonzero entries, satisfying \ref{triangle} and \ref{square}.  Let $\Mm^{\sigma}$ be the $N$-scheme obtained after a permutation $\sigma$ of the indices (replacing $m_{\sigma(i) \sigma(j)}$ by $-m_{\sigma(j) \sigma(i)}$ if $\sigma(i)>\sigma(j)$). Then \\
(1) $\toz(\Mm;p)<p$ for all $p<N$ if and only if $\toz(\Mm^{\sigma};p)<p$ for all $p<N$. \\
(2) In that case replacing the definition (\ref{toz}) with the following more restrictive one, the main theorem still holds 
\begin{equation*}
\toz(\Mm;g_l,j)= \left\{ \begin{array}{ll} 
(\nu_{lj}-\nu_{l}+1)^{-1} & \nu_{l1j}=\ldots=\nu_{l(j-1)j}\leq \nu_{l}, \\
0 & \mbox{otherwise}.
\end{array} \right.  \end{equation*}
\end{cor}

\section{Examples and Computations} \label{section:examp}

Now that we know how to detect when a given $\Mm$ can be realized on a torus, we  produce examples. 

\subsection{Intersection numbers bounded from above.} \label{section:agol}
In \cite{Ag1} the following question is considered. Suppose a positive integer $d$ is given. Let $\Delta_{\leq} (d)$ be the maximal number of curves which lie on a torus with the condition that absolute values of their pairwise intersection numbers are bounded above by $d$. 
Then Agol proves a lower bound and conjectures that the number $P(d)=1+$~(next prime greater than $d$)  is an upper bound for $\Delta_{\leq} (d)$. In cases where these two bounds differ, finer lower and upper bounds are proposed via computer aid. 

Theorem~\ref{torus} can be utilized to attack this problem. In the simplest case $d=7$ discussed in \cite{Ag1}, we can also confirm  via a short C code that $\Delta_{\leq}(d)=10$. The next case is when $d=13$ and our naive code cannot handle such big data in a reasonable time.

\subsection{A 2-endemic 4-scheme} \label{section:endem}
Consider a 4-scheme  $\Mm$ with all entries odd. $g(\Mm)>1$ since it cannot satisfy \ref{square}. 
Fix two prime numbers $p,q>2$.  Assume that all entries of $\Mm$ are divisible by $p$ and $q$ except $p\nmid m_{12}$ and $q\nmid m_{34}$.
We want to express $\Mm$ as the sum of two 4-schemes $\Mm'$ and $\Mm''$ so that both satisfy \ref{circledast2} and \ref{square}. Let $E'$ and $E''$ be the number of their even entries respectively. Note that their corresponding entries are of different parity so that $E'+E''=6$. So $E'$ cannot be 0 or 6. Moreover \ref{circledast2} prohibits the other values except the cases when some entries are 0.

So in all these cases let us investigate the possibility of having zeros in $\Mm'$ and $\Mm''$.
First suppose $E'=5$ and $\Mm'=(0;1,0;0,0,0)\mod 2$. To satisfy \ref{circledast2} we must have at least one of 
$m'_{12},m'_{23}$ zero. If just $m'_{12}=0$ then $\gamma'_1\sim \pm \gamma'_2$ (i.e. isotopic) so that $|m'_{13}|=|m'_{23}|$ which is impossible.
So both $m'_{12}=m'_{23}=0$. Similar discussion for the triple  $1,3,4$ forces $m'_{14}=m'_{34}=0$ too. But in this case $m'_{24}=0$ too (and both $\gamma'_2$, $\gamma'_4$ are  empty). In that case \ref{circledast2} is not satisfied in $\Mm''$ 
for indices $2,3,4$ (by assumption on divisibility by $p$).
Same discussion goes through for any other value for $\Mm'$ except $(0;0,0;0,0,1)\mod 2$ in case of which \ref{circledast2}
is not satisfied in $\Mm''$ either for indices $1,2,3$  (by assumption on divisibility by $q$).

Now suppose for $E'=4$, without loss of generality,  $\Mm'=(0;0,1;1,1,1)\mod 2$. To satisfy \ref{circledast2} in $\Mm'$ we must have at least one of 
$m'_{12},m'_{13}$ zero. As before both are forced to be zero. This in turn requires either $m'_{23}=0$ ($\gamma'_1\sim \pm \gamma'_2\sim \pm \gamma'_3$) 
or $m'_{14}=0$ ($\gamma'_1$ empty), which are impossible.

Finally suppose $E'=3$. Considering \ref{circledast2} and \ref{square} the only possible setup is when, without loss of generality,  
either $\Mm'=(0;0,1;0,1,1)\mod 2$ or $\Mm'=(1;1,1;0,0,0)\mod 2$. (We have two cases because of the asymmetry caused by the assumption on divisibility by $p$ and $q$.) 
Either is possible only if all even entries are zero. In the former case $\gamma'_1$ becomes empty.
However in that case all $m''_{23},m''_{24},m''_{34}$ are to be divisible by $q$ because of the initial assumption on $m_{12},m_{13},m_{14}$.
Then $q|m'_{23},q|m'_{24}$ while $q\nmid m'_{34}$. Thus $\Mm'$ fails to satisfy \ref{circledast2}.
(Note that the argument works regardless of $m''_{23},m''_{24}$ or $m''_{34}$ being zero.)

In the latter case $\gamma'_4$ becomes empty. 
However in that case $m''_{12},m''_{13},m''_{23}$ are to be divisible by $p$ because of the initial assumption on $m_{14},m_{24},m_{34}$.
This implies $p|m'_{12},p|m'_{13},p\nmid m'_{23}$. Thus $\Mm'$ fails to satisfy \ref{circledast2}. \\

\begin{figure}[h]
\centering
\includegraphics[width=0.3\columnwidth]{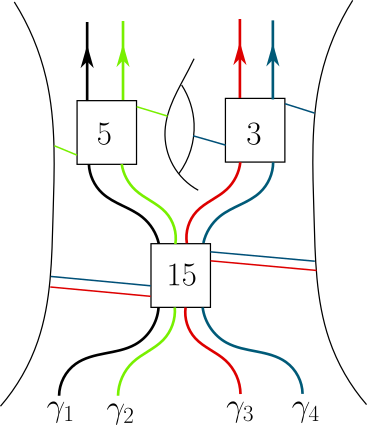}
\caption{The 4-scheme $\Mm=\left(5;15,15;15,15,3 \right)$ on $\Sigma_2$; top and bottom are identified. The boxes are as described in Figure~\ref{babayuzey}.}
\label{2demik}
\end{figure}

\noindent {\em Proof of Theorem~\ref{triples}(b).}
For prime numbers $p,q>2$ consider  the 4-scheme 
$\Mm=\left(q;pq,pq;pq,pq,p \right)$ which does not satisfy neither \ref{circledast2} nor \ref{square}. 
So we have $g(\Mm)>1$. However it is immediate to observe that  $g(\Mm)=2$ (see Figure~\ref{2demik}).
The discussion above shows that $\Mm$ is in fact a 2-endemic 4-scheme. \hfill $\Box$

\subsection{Various examples} \label{section:varexamp}

Here are some examples in each of which $\Mm$ satisfies \ref{triangle} and  \ref{square}.
We display the schemes in matrix form to enhance legibility. Recall that for each prime divisor $p$ of $g_{123}$, a column has a nonzero contribution to $\toz(\Mm;p)$ if and only if each entry in that column has the same $p$-valuation.

We start with a 4-scheme which is realizable on a torus:

\begin{tabular}{cc}
\multirow{3}{*}{
$\Mm=\left(m_{ij}\right)_{1\leq i<j\leq 4}:\;\;${\footnotesize $\begin{matrix}
1&1&2\\
&1&1\\
&&-1
\end{matrix}$}} & \\
& $\toz(\Mm;2)=0, \toz(\Mm;3)=0$. \\
&
\end{tabular} \\

\noindent Multiplying all entries with 3 makes $\toz(\Mm;3)=1+1+1$, failing \ref{circledast}.

\begin{tabular}{cc}
\multirow{3}{*}{
{\footnotesize $\begin{matrix}
3&3&6\\
&3&3\\
&&-3
\end{matrix}$}} & \\
& $\toz(\Mm;2)=0, \toz(\Mm;3)=3$. \\
&
\end{tabular} \\


\noindent Multiplying first 3 entries of this scheme with 3 makes $\toz(\Mm;3)=1+1+1/3$. So the scheme can be realized on a torus.

\begin{tabular}{cc}
\multirow{3}{*}{
{\footnotesize $\begin{matrix}
9&9&6\\
&9&3\\
&&-3
\end{matrix}$}} & \\
& $\toz(\Mm;2)=0, \toz(\Mm;3)=2+\frac{1}{3}$. \\
&
\end{tabular} \\


\noindent The following 6-scheme satisfies \ref{circledast} and hence can be realized on a torus.

\begin{tabular}{cc}
\multirow{5}{*}{
{\footnotesize $\begin{matrix}
3&3&1&-1&1\\
&6&4&2&1\\
&&2&4&-1 \\
&&&2&-1 \\
&&&&-1 \\
\end{matrix}$}} & \\
& \\
& $\toz(\Mm;2)=0, \toz(\Mm;3)=1, \toz(\Mm;5)=0$. \\
& \\
& \\
\end{tabular} \\


\noindent Multiplying all entries with 3 makes $\toz(\Mm;3)=1+1+\frac{1}{3}+\frac{1}{3}+\frac{1}{3}$. This fails \ref{circledast}.

\begin{tabular}{cc}
\multirow{5}{*}{
{\footnotesize $\begin{matrix}
9&9&3&-3&3\\
&18&12&6&3\\
&&6&12&-3 \\
&&&6&-3 \\
&&&&-3 \\
\end{matrix}$}} & \\
& \\
& $\toz(\Mm;2)=0, \toz(\Mm;3)=3, \toz(\Mm;5)=0$. \\
& \\
& \\
\end{tabular} \\


\noindent Instead, multiplying all entries with 5 makes $\toz(\Mm;5)=1+1+1+1+1$, failing \ref{circledast}.

\begin{tabular}{cc}
\multirow{5}{*}{
{\footnotesize $\begin{matrix}
15&15&5&-5&5\\
&30&20&10&5\\
&&10&20&-5 \\
&&&10&-5 \\
&&&&-5 \\
\end{matrix}$}} & \\
& \\
& $\toz(\Mm;2)=0, \toz(\Mm;3)=1, \toz(\Mm;5)=5$. \\
& \\
& \\
\end{tabular} \\


\noindent Here is another 6-scheme with a factor 5 in each entry. Different than the previous example, the existence of an extra 5 factor of $m_{23}$ makes $\toz(\Mm;5)=1+1+1+1$. So \ref{circledast} holds and the scheme can be realized on a torus.

\begin{tabular}{cc}
\multirow{5}{*}{
{\footnotesize $\begin{matrix}
15&20&5&-15&-10\\
&\mathbf{25}&10&15&-5\\
&&5&5&10 \\
&&&5&5 \\
&&&&-5 \\
\end{matrix}$}} & \\
& \\
& $\toz(\Mm;2)=0, \toz(\Mm;3)=0, \toz(\Mm;5)=4$. \\
& \\
& \\
\end{tabular} \\


\noindent Replacing one factor 5 with 3 in each entry makes $\toz(\Mm;3)=1+1+1+1+1$, failing \ref{circledast}.

\begin{tabular}{cc}
\multirow{5}{*}{
{\footnotesize $\begin{matrix}
9&12&3&-9&-6\\
&15&6&9&-3\\
&&3&3&6 \\
&&&3&3 \\
&&&&-3 \\
\end{matrix}$}} & \\
& \\
& $\toz(\Mm;2)=0, \toz(\Mm;3)=5, \toz(\Mm;5)=0$. \\
& \\
& \\
\end{tabular} \\



\end{document}